\documentclass[11pt]{amsart}
\usepackage[square,compress,comma, numbers]{natbib}
\usepackage[colorlinks=true, citecolor=blue, linkcolor=blue]{hyperref}
\usepackage{amsfonts}
\allowdisplaybreaks[4]
\usepackage{amssymb}
\usepackage{amsmath}
\usepackage{color}
\newcommand{\abs}[1]{\left\lvert #1 \right\rvert}

\def\E#1{\mathbb{E}\left \{#1 \right\}}

\definecolor{c20}{rgb}{0.,0.7,0.}
\definecolor{c30}{rgb}{0.,0.,1.}
\definecolor{c40}{rgb}{1,0.1,0.7}
\definecolor{c50}{rgb}{1,0,0}
\definecolor{c60}{rgb}{1,0.9,0.1}
\definecolor{c70}{rgb}{0.50,1.00,0.00}

\numberwithin{equation}{section}
\newtheorem{theo}{Theorem}[section]
\newtheorem{sat}[theo]{Proposition}
\newtheorem{de}[theo]{Definition}
\newtheorem{lem}{Lemma}[section]

\newtheorem{example}[theo]{Example}
\newtheorem{korr}[theo]{Corollary}
\newtheorem{remark}[theo]{Remark}
\newtheorem{remarks}[theo]{Remarks}
\numberwithin{equation}{section}

\newcommand{\prooftheo}[1]{ \textsc{Proof of Theorem} \ref{#1} }

\newcommand{\prooflem}[1]{\textsc{Proof of Lemma} \ref{#1}}

\newcommand{\pk}[1]{\mathbb{P} \left\{ #1 \right\} }

\newcommand{\QED}{\hfill $\Box$}
\newcommand{\COM}[1]{}
\def\IF{\infty}

\newcommand{\BQN}{\begin{eqnarray}}
\newcommand{\EQN}{\end{eqnarray}}
\newcommand{\BQNY}{\begin{eqnarray*}}
\newcommand{\EQNY}{\end{eqnarray*}}

\def\polhk#1{\setbox0=\hbox{#1}{\ooalign{\hidewidth
\lower1.5ex\hbox{`}\hidewidth\crcr\unhbox0}}}

\def\cL#1{\textcolor{c50}{#1}}
\def\cL#1{#1}

\newcommand{\ve}{\varepsilon}

\def\vp{\varepsilon}

\def\rw{\rightarrow}

\def\IF{\infty}

\def\Cov{\mathrm{Cov}}

\date{}

\def\LT{\left}
\def\RT{\right}

\def\ooo{(1+o(1))}

\def\rw{\rightarrow}

\def\Piter{\mathcal{P}}

\def\TT{\mathcal{T}}

\def\vn{\varepsilon}
\def\Var{\text{Var}}

\def\FF{\widetilde{\mathcal{H}}}

\newcommand{\BS}{\begin{sat}}
\newcommand{\ES}{\end{sat}}
\newcommand{\BT}{\begin{theo}}
\newcommand{\ET}{\end{theo}}
\newcommand{\BK}{\begin{korr}}
\newcommand{\EK}{\end{korr}}

\newcommand{\BD}{\begin{de}}
\newcommand{\ED}{\end{de}}
\newcommand{\BIT}{\begin{itemize}}
\newcommand{\EIT}{\end{itemize}}
\newcommand{\BDI}{\begin{description}}
\newcommand{\EDI}{\end{description}}
\newcommand{\BRM}{\begin{remarks}}
\newcommand{\ERM}{\end{remarks}}
\newcommand{\BEL}{\begin{lem}}
\newcommand{\EEL}{\end{lem}}

\def\TT{\mathcal{T}}
\def\TT{\R }

\def\rw{\rightarrow}
\def\LT{\left}
\def\RT{\right}

\def\Var{\text{Var}}
\def\vn{\varepsilon}
\def\Cov{\mathrm{Cov}}

\def\LT{\left}
\def\RT{\right}

\def\ooo{(1+o(1))}

\def\rw{\rightarrow}

\def\Piter{\mathcal{P}}

\def\TT{\mathcal{T}}

\def\vn{\varepsilon}
\def\Var{\text{Var}}

\def\LT{\left}
\def\RT{\right}

\def\ooo{(1+o(1))}

\def\rw{\rightarrow}

\def\Piter{\mathcal{P}}

\def\TT{\mathcal{T}}

\def\vn{\varepsilon}
\def\Var{\text{Var}}

\def\FF{\widetilde{\mathcal{H}}}

\def\cL#1{\textcolor{c50}{#1}}
\def\cL#1{#1}

\begin{document}
\title{Parisian Ruin Probability of an Integrated Gaussian Risk Model}
\author{Xiaofan Peng}
\address{Xiaofan Peng, School of Mathematical Sciences, University of Electronic Science and Technology of China, Chengdu 610054, China}
\email{xfpengnk@126.com}
\author{Li Luo}
\address{Li Luo, School of Mathematical Sciences, Nankai University, Tianjin 300071, P.R. China and Department of Actuarial Science, University of Lausanne, UNIL-Dorigny, 1015 Lausanne, Switzerland}
\email{luol14@mail.nankai.edu.cn}
\bigskip
\maketitle
{\bf Abstract:} In this paper we investigate the Parisian ruin probability for an integrated Gaussian process. Under certain assumptions, we find the Parisian ruin probability and the classical ruin probability are on the log-scale asymptotically the same. Moreover, for any small interval required by the risk process staying below level zero, the Parisian ruin probability and the classical one are the same also in the premise asymptotic behavior. Furthermore, we derive an approximation of the conditional ruin time. 
\def\TTT{\mathcal{T}}
\def\TT{\mathcal{T}}
\def\Z{\mathbb{Z}}
\def\HWD{\mathcal{H}_\alpha^\delta}
\def\phd{\mathcal{P}_{\alpha,\delta}^h}

\section{Introduction}
 Gaussian risk processes have been investigated in numerous research paper. With motivation from \cite{rolski2009stochastic}  the risk
 reserve process of an insurance company can be modelled by a
 stochastic process $\{R_u(t),t\ge 0\}$ given as
 \BQNY\label{U0}
 R_u(t)=u+ct-\int_0^tZ(s)\, ds,\ \ \ t\ge0,
 \EQNY where $u\ge0$ is the initial reserve, $c>0$ is the rate of
 premium received by the insurance company  and $\{Z(t),t\ge0\}$ is a centered Gaussian process with almost surely continuous sample paths.
 Commonly  the process $\{Z(t),t\ge0\}$ is {referred to as} the loss rate of the insurance company.
In order to take into account the time-value of money, in this contribution for a given
real-valued measurable function $\delta(\cdot)$  we shall consider the more general risk process
\BQN\label{eq:model}
R_u(t)=u+c\int^{t}_{0}e^{-\delta(s)}d s-\int_{0}^{t}e^{-\delta(s)}Z(s)d s, \quad t\geq0.
\EQN
An important quantity of interest for such a risk process is the calculation of the ruin probability over the finite time-horizon $[0,S]$
\begin{align*}
\psi_S(u)&:=\pk{\inf_{t\in[0,S]}R_u(t)<0}\\
&=\pk{\sup_{t\in[0,S]}\biggl(\int_0^t e^{- \delta(s)} Z(s)\, ds-c\int_0^t e^{-\delta(s)}ds\biggr)>u}.
\end{align*}
Since it is not possible to calculate $\psi_S(u)$ for any fixed $u$ explicitly, one resorts to asymptotic theory analysing the ruin probability as the initial reserve $u$ becomes large. The recent contribution \cite{dbicki2015gaussian} derived the exact tail asymptotics of $\psi_S(u)$ as $u\to \IF$ under some restrictions on $Z$. Moreover, therein an approximation of the ruin time as $u\to \IF$ is derived.

A more general concept than the classical ruin probability is the Parisian ruin probability, which in our context is defined for some given $T_u>0$
as
\begin{align}\label{Intruin1}
\mathcal{P}_S(u,T_{u})&:=\pk{\inf_{t\in[0,S]}\sup_{s\in[t,t+T_{u}]}R_u(s)<0}.
\end{align}
Clearly, in the particular case that $T_u=0$ we have that the Parisian ruin probability equals the classical ruin probability.
Parisian ruin and Parisian ruin time have been recently discussed for self-similar Gaussian processes in \cite{debicki2015parisian,debicki2016parisian}, whereas the recent publication \cite{BaiL} discusses the classical Brownian motion risk model. With motivation from these recent contributions we shall analyse the
asymptotic behavior of $\mathcal{P}_S(u,T_{u})$ as $u\to \IF$. Our objectives are two-fold. First we are interested in a large deviation type result for  the Parisian ruin probability. Specifically, under two weak restrictions on $Z$ and the inflation/deflation rate function $\delta(\cdot)$ we shall show in our first result that
\BQNY
\lim_{u\to\infty}\frac{\log(\Piter_S(u,\cL{T_u}))}{u^2}=-\frac{1}{2\sigma^2(S)},
\EQNY
  where $\sigma^2(t)$ is given by
  \BQN \label{sigm}
  \sigma^2(t):=\Var(R_{u}(t))=2\int_0^t\int_{0}^{v}e^{-\delta (w)-\delta (v)}\Cov(Z(w),Z(v))dwdv.
  \EQN

 The above result is shown for quite general $T_u$, in particular it holds for $T_u=0$. Hence,
  the Parisian ruin probability and the classical ruin probability are on the log-scale asymptotically the same, i.e.,
  \BQNY \label{fristRes}
\log(\Piter_S(u,\cL{T_u})) \sim \log \psi_S(u), \quad u\to \IF,
\EQNY
where $\sim$ stands for asymptotic equivalence when $u\to \IF$.  Then we show in our main result, for any small interval that required by the risk reserve process staying below level zero, it is possible to derive the exact asymptotic of the Parisian ruin probability as $u\to \IF$. Such a result reveals, that the asymptotic of Parisian ruin probability and the classical one are the same also in the precise asymptotic behavior.

Brief organization of the rest of the paper: Section \ref{section-re} presents our main results where additionally to the large deviation and the precise asymptotic of the Parisian ruin probability. Furthermore, we obtain an approximation of the Parisian ruin time. All the proofs are relegated to Section \ref{section-proof} which concludes this contribution.

\section{Results}\label{section-re}
In this section we shall present two results. The first one gives the large deviation asymptotic and the precise asymptotic of the Parisian ruin probability.
The second result is concerned with the Parisian ruin time. In theoretical investigations, the analysis of ruin time is of interests since it gives more information on how and when the ruin occurs. Our results are derived under the following conditions on the risk reserve process $R_u(t)$.

\textbf{A1}. The claim rate process $Z(t)$ is a centered, non-degenerate Gaussian process with continuous sample path and nonnegative covariance,
 i.e., $\Cov(Z(s),Z(t))\geq0$ for any $s,t\geq0$.\\
\textbf{A2}. The inflation/deflation rate function $\delta(\cdot)$ is locally bounded.

Clearly, condition \textbf{A2} is always met in practical applications. Condition \textbf{A1} is a weak one, it is satisfied by many Gaussian processes, for instance, Ornstein-Ohlenbeck process, Slepian process and the fractional Brownian motion with Hurst index $H \in (0, 1)$. Note that $H=1/2$ corresponds to the Brownian motion.

Next, we give two asymptotic results which are based on the characteristics of pre-specified time $T_u$.

\BT\label{pan1}
Let $\{R_u(t),t\ge0\}$ be the reserve process evolving as (\ref{eq:model}) and satisfy assumptions \textbf{A1-A2}. Define
$\sigma^2(\cdot)$ by \eqref{sigm} and set $\tilde\delta (t):=\int_{0}^{t}e^{-\delta(s)}ds$. For any bounded delayed time $T_{u}\geq0$ we have
\BQNY
\lim_{u\to\infty}\frac{\log(\Piter_S(u,\cL{T_u}))}{u^2}=\lim_{u\to\infty}\frac{\log(\psi_S(u))}{u^2}=-\frac{1}{2\sigma^2(S)}.
\EQNY
Furthermore, for any small interval $T_{u}\to0$ as $u\to\infty$, then
\BQNY\label{main}
\mathcal{P}_S(u,T_u)\ \sim \ \psi_S(u)\ \sim  \ \pk{ \mathcal{N}>(u+ c \tilde\delta (S))/\sigma(S)}
\EQNY
holds as $u\rightarrow \infty$, with $\mathcal{N}$ a $N(0,1)$ random variable.
\ET

Another quantity of interest is the conditional distribution of the ruin time for the surplus process $R_u(s)$. We define the ruin time as
\BQN\label{eq:eta2}
\tau(u):=\inf\{t\ge T_u: t-\kappa_{t,u}\ge T_u, R_u(t)<0\},\ \ \ \text{with} \ \kappa_{t,u}=\sup\{s\in[0,t]: R_u(s)\ge0\}.
\EQN
\BT\label{pan3}
Under the conditions of Theorem \ref{pan1}, for small interval $T_u$, we have
\BQN\label{tau}
\lim_{u\to\infty}
\pk{u^2(S+T_{u}-\tau(u))\leq x \Bigl|\tau(u)<S+T_{u}}
=
1 - \exp\left(     -  \frac{\sigma'(S)}{\sigma^3(S)}x    \right),\ \ x\ge0.
\EQN
\ET
Below we shall present two illustrating examples. It is worth noting that all the examples given in \cite{dbicki2015gaussian} are adapted to our model.
\begin{example}
Let $\{Z(t), t\geq0\}$ be a standard fractional Brownian motion with Hurst index $H\in(0,1)$, i.e., it is a centered Gaussian process with a.s. continuous sample paths and covariance function $\Cov(Z(t),Z(s))=\frac{1}{2}(\abs{t}^{2H}+\abs{s}^{2H}-\abs{t-s}^{2H})$ . If $\delta(t)=t, t\geq0$, then
\begin{align*}
\widetilde{\delta}(t)&=1-e^{-t}, \\
\sigma^2(t)&=2\int_{0}^{t}\int_{0}^{x}e^{-x-y}(x^{2H}+y^{2H}-(x-y)^{2H})dydx\\
&=\Gamma(2H+1,t)(1-2e^{-t})+e^{-2t}\Gamma^*(2H+1,t),
\end{align*}
where $\Gamma(a,t)=\int_0^{t}x^{a-1}e^{-x}dx$ and $\Gamma^*(a,t)=\int_0^{t}x^{a-1}e^{x}dx$. It is easy to see conditions $\mathbf{A1}$ and $\mathbf{A2}$ are naturally satisfied. Consequently, Theorem \ref{pan1} implies, for $T_u\to0$ as $u\to\IF$,
\begin{align*}
&\mathcal{P}_S(u,T_u)=\frac{1}{u}\sqrt{\frac{\Gamma(2H+1,S)(1-2e^{-S})+e^{-2S}\Gamma^*(2H+1,S)}{2\pi}}\\
&\qquad\qquad\quad\times \exp\LT(-\frac{(u+c(1-e^{-S}))^2}{2\Gamma(2H+1,S)(1-2e^{-S})+2e^{-2S}\Gamma^*(2H+1,S)}\RT)(1+o(1)).
\end{align*}
Furthermore, according to Theorem \ref{pan3} the convergence in \eqref{tau} holds with
\BQNY
\frac{\sigma'(S)}{\sigma^3(S)} = \frac{e^{-S}\LT(S^{2H}+\Gamma(2H+1,S)\RT) -e^{-2S}\LT( S^{2H}+\Gamma^*(2H+1,S)\RT)}{(\Gamma(2H+1,S)(1-2e^{-S})+e^{-2S}\Gamma^*(2H+1,S))^2}.
\EQNY
\end{example}
\begin{example}
Let $\{Z(t)=\frac{B(t)}{\sqrt{t}}, t\geq0\}$ be a scaling Brownian motion with $B$ a standard Brownian motion. If further $\delta(t)=t, t\geq0$, then by Taylor formula
\begin{align*}
\sigma^2(t)&=2\int_{0}^{t}\int_{0}^{x}e^{-x-y}\sqrt{y/x}\ dydx\\
&=2\int_{0}^{1}\frac{\sqrt{z}(1-e^{-t(1+z)})}{(1+z)^2}dz - 2\int_{0}^{1}\frac{t\sqrt{z}e^{-t(1+z)}}{1+z}dz \\
&=\frac{2}{3}t^2+\sum_{k=3}^{\infty}(-1)^{k}\frac{2(k-1)}{k!}t^{k}\int_{0}^{1}(1+z)^{k-2}\sqrt{z}dz.
\end{align*}
Applying Theorem \ref{pan1} once again, we obtain
\begin{align*}
&\mathcal{P}_S(u,T_u)=\frac{1}{u}\sqrt{\frac{\frac{2}{3}S^2+\sum_{k=3}^{\infty}(-1)^{k}\frac{2(k-1)}{k!}S^{k}\int_{0}^{1}(1+z)^{k-2}\sqrt{z}dz}{2\pi}}\\
&\qquad\qquad\quad\times \exp\LT(-\frac{(u+c(1-e^{-S}))^2}{\frac{4}{3}S^2+\sum_{k=3}^{\infty}(-1)^{k}\frac{4(k-1)}{k!}S^{k}\int_{0}^{1}(1+z)^{k-2}\sqrt{z}dz}\RT)(1+o(1)),
\end{align*}
for $T_u\to0$ as $u\to\IF$. Finally, by Theorem \ref{pan3} the convergence in \eqref{tau} holds with
\BQNY
\frac{\sigma'(S)}{\sigma^3(S)} = \frac{S^{-\frac{1}{2}}e^{-S}\Gamma(\frac{3}{2},S)}
{(\frac{2}{3}S^2+\sum_{k=3}^{\infty}(-1)^{k}\frac{2(k-1)}{k!}S^{k}\int_{0}^{1}(1+z)^{k-2}\sqrt{z}dz)^2}.
\EQNY
\end{example}

\section{Proofs}\label{section-proof}

Before the demonstration, for notational simplicity, we define
\begin{align*}
 Y(t)&:=\int_{0}^{t}e^{-\delta (s)}Z(s)ds,  \qquad\ \   R(s,t):=\Cov(Z(s),Z(t))\\
g_u(t)&:= \frac{u+c\tilde\delta (t)}{\sigma(t)},  \qquad\qquad\qquad X_u(t):=\frac{Y(t)}{\sigma(t)}\frac{g_u(S)}{g_u(t)},\\
\sigma_{X_u}^2(t)&:=\Var(X_u(t)),\qquad\qquad  r_{X_u}(s,t):=\mathbb{C}orr\LT(X_u(s),X_u(t)\RT).
\end{align*}
Then, we can reformulate (\ref{Intruin1}) as
\BQN\label{psu-ano-form}
\mathcal{P}_S(u,T_u)=\pk{\sup_{t\in[0,S]}\inf_{s\in[t,t+T_u]}X_u(s) > g_u(S)}.
\EQN
\prooftheo{pan1}
From assumption \textbf{A1}, we know in fact $Z(t)$ is continuous in the mean squared sense. This means that the covariance function $R(s,t)$ is
a bivariate continuous function and is strictly positive for $|t-s|$ sufficiently small due to the non-degeneracy. Therefore,
\BQN\label{eq:dif-vari}
\frac{\partial\sigma^{2}(t)}{\partial t}=2\int_{0}^{t}e^{-\delta(s)-\delta(t)}R(s,t)ds >0
\EQN
 and for $0<s\leq t$
\BQNY\label{eq:compa-vari}
\Cov(Y(t),Y(s))\geq\sigma^{2}(s).
\EQNY
Then by the Slepian Lemma (cf. \cite{piterbarg1996asymptotic})
\BQNY
\mathcal{P}_S(u,T_u)
\le\pk{\sup_{t\in[0,S]} Y(t)>u} \le \pk{\sup_{t\in[0,S]} B(\sigma^{2}(t))>u}=2\Psi\left(\frac{u}{\sigma(S)}\right),
\EQNY
where $B(\cdot)$ is a standard Brownian motion and $\Psi(\cdot)$ denotes the tail distribution of a standard normal random variable. Moreover, appealing to Theorem 2.1 in \cite{dbicki2015gaussian}, we get a sharper upper bound
\BQN\label{upp-bound}
\mathcal{P}_S(u,T_u) \leq \psi_S(u) = \Psi(g_u(S))\ooo, \qquad \textmd{as}\ u\rw\IF.
\EQN
For the lower bound, put $T=\sup_{u>0}T_{u}$ and $m=\sup_{t\in[0,S+T]}e^{-\delta(t)}$. Due to our assumptions, $T$ and $m$ are both finite. Then drawing on similar arguments as used in the proof of Theorem 3.1 in \cite{debicki2016parisian} (only replacing therein $c$ by $cm$ and $\rho_{S}$ by $(\inf_{t\in[S,S+T]}\partial\sigma^2(t)/\partial t)^{-1}$, note that the latter is finite and the concavity in the aforementioned paper is not necessary) we get
\BQNY
\mathcal{P}_S(u,T_u)
\ge\pk{\inf_{t\in[S,S+T]} Y(t)-cmt>u}\geq C\frac{\sigma^2(S)}{u} \Psi\LT(\frac{ u+cmS}{\sigma(S)}\RT)\ooo
\EQNY
for some positive constant $C$, as $u\to\IF$.\\
The first claim follows straightforwardly from combination of the above inequalities concerning $\mathcal{P}_S(u,T_u)$.\\
Next, we derive a lower bound of $\mathcal{P}_S(u,T_u)$ for $T_{u}\to0$ as $u\to\infty$ by considering three separate cases.
%

\textbf{Case I:} $T_{u}=o(\frac{1}{u})$\label{case-i}

 First, differentiating $\sigma_{X_u}(s)$ yields
\BQN\label{diff-xu}
\sigma'_{X_u}(s)=\frac{\sigma'(s)}{\sigma(S)}\frac{u+c\tilde\delta (S)}{u+c\tilde\delta (s)}-\frac{ce^{-\delta(s)}\sigma(s)(u+c\tilde\delta (S))}{(u+c\tilde\delta (s))^2\sigma(S)},
\EQN
which together with (\ref{eq:dif-vari}) implies $\sigma'_{X_u}>0$ for sufficiently large $u$. Consequently, $\sigma_{X_u}(s)>1$ for all $s>S$. Secondly, given arbitrary $\theta>S$, then for any $s,t\in[S,\theta]$
\begin{align}\nonumber
1-r_{X_u}(s,t)&\leq\frac{\Var(Y(t)-Y(s))}{2\sigma(s)\sigma(t)}\\\nonumber
&=\frac{\int_s^t\int_s^t  \cL{R(w,v)}e^{-\delta(w)-\delta(v)} dwdv}
     {2\sigma(s)\sigma(t)}\\ \label{upp-bound-corr}
     &\leq C(t-s)^2,
\end{align}
where $C=\max_{w,v\in[S,\theta]}R(w,v)e^{-\delta(w)-\delta(v)}/(2\sigma^{2}(S))$. Next, for any $\varepsilon>0$, put $C_\vp=C(1+\vp)$ and define a centered Gaussian process $\{\xi_\ve(t),t\ge0\}$ with covariance function
$\Cov(\xi_\ve(t),\xi_\ve(s))=e^{-C_\ve(t-s)^2}$. Lastly, in view of Slepian Lemma, for $T_{u}=o(u^{-1})$ and any sufficiently small $\varepsilon_{1}\in(0,1)$
\begin{align}
\mathcal{P}_S(u,T_u)&\geq\pk{\inf_{s\in[S,S+T_{u}]} X_u(s)\ >\ g_u(S)}\nonumber\\
&\geq\pk{\inf_{s\in[S,S+\varepsilon_{1}u^{-1}]} \frac{X_u(s)}{\sigma_{X_u}(s)}\ >\ g_u(S)}\nonumber\\
&\geq\pk{\inf_{s\in[0,\varepsilon_{1}]}\xi_\ve(su^{-1})>g_u(S)}\nonumber\\
&=\FF_2(\hat a\varepsilon_{1}) \Psi\LT(g_u(S) \RT)\ooo\label{small-s-Tu}
\end{align}
as $u$ sufficiently large, where the last equality follows from Lemma 5.1 in \cite{debicki2016parisian}, $\widehat{a}=\sqrt{C_\vp}/\sigma(S)$ and $\FF_2(\cdot)$ is described by the following generalized Pickands constant
\BQNY\label{pickand-cons}
\FF_\alpha(T)=\E{\exp\LT(\inf_{ s\in[0,T]}\LT(\sqrt{2}B_{\alpha}(s)-s^{\alpha}\RT)\RT)}\in(0,\IF),\ \ T\ge0,
\EQNY
with $B_{\alpha}(\cdot)$ a standard fractional Brownian motion with Hurst index $\alpha/2\in(0,1]$.
Therefore, letting $\vn_1\to 0$ in \eqref{small-s-Tu} yields the lower bound
\BQN\label{lower-bound}
\mathcal{P}_S(u,T_u)\geq\Psi\LT(g_u(S) \RT)\ooo, \quad u\to\IF.
\EQN

\textbf{Case II:} $T_{u}=O(\frac{1}{u})$

Without loss of generality, suppose $T_{u}=T/u$ for some positive constant $T$. Then, similar to the case I, for any constant $Q>0$ and sufficiently large $u\ (>2Q\sigma'(S)/\sigma(S))$
\begin{align}
\mathcal{P}_S(u,T_u)&\geq\pk{\inf_{s\in[S,S+T_{u}]} X_u(s)\ >\ g_u(S)}\nonumber\\
&\geq\pk{\inf_{s\in[S,S+Tu^{-1}]} \frac{X_u(s)}{\sigma_{X_u}(s)}\LT(1+\frac{Q}{u}(s-S)\RT)\ >\ g_u(S)}\nonumber\\
&\geq\pk{\inf_{s\in[0,T]}\xi_\ve\LT(\frac{s}{u}\RT)\LT(1+\frac{Qs}{u^2}\RT)>g_u(S)}\nonumber\\
&=\FF_2^{Q/\sqrt{C_\vp}}(\hat{a}T) \Psi\LT(g_u(S)\RT)\ooo\label{small-Tu},
\end{align}
where $\FF_2^{Q}(\cdot)$ is described by the following generalized Piterbarg constant
\BQNY\label{pickand-cons}
\FF_\alpha^Q(T)=\E{\exp\LT(\inf_{ s\in[0,T]}\LT(\sqrt{2}B_{\alpha}(s)-s^{\alpha}+Qs\RT)\RT)}\in(0,\IF),\ \ T\ge0,\ \alpha\in(0,2].
\EQNY
Letting $Q\to\infty$ in \eqref{small-Tu} gives the same lower bound as (\ref{lower-bound}).

\textbf{Case III:} $\frac{1}{u}=o(T_{u})$ and $T_u\to0$ as $u\to\infty$

From the proof above we see that the small constant $\vp$ plays an insignificant role. Hence, for the sake of notational convenience, we use
 $C$ and $\xi(\cdot)$ instead of $C_\vp$ and $\xi_\ve(\cdot)$,
and put $\hat{\xi}(s)=\xi(s)\sigma_{X_u}(s+S)$ for $s\geq0$. Then, using Slepian Lemma again yields
\begin{align}
\mathcal{P}_S(u,T_u)&\geq\pk{\inf_{s\in[S,S+T_{u}]} X_u(s)\ >\ g_u(S)}\nonumber\\
&\geq\pk{\inf_{s\in[0,T_u]}\hat{\xi}(s)>g_u(S)}\nonumber\\\label{primai-prob}
&=\pk{\hat{\xi}(0)>g_u(S),\ \hat{\xi}(T_u)>g_u(S)}\\\label{negli-prob}
&\qquad-\ \pk{\hat{\xi}(0)>g_u(S),\ \exists s\in(0,T_u)\ s.t.\ \hat{\xi}(s)\leq g_u(S),\ \hat{\xi}(T_u)>g_u(S)}.
\end{align}
In the following, we first show (\ref{primai-prob}) is asymptotically equivalent to $\Psi\LT(g_u(S)\RT)$ as $u\to\infty$, and then appeal to the method of moments (see \cite{piterbarg1996asymptotic}) to show that (\ref{negli-prob}) is negligible with respect to the former probability. Note, for a bivariate normal random variable,
\begin{align}
&\pk{\hat{\xi}(0)>g_u(S),\ \hat{\xi}(T_u)>g_u(S)}\nonumber\\\nonumber
&=\mathbb{P}\Big\{\xi(0)>g_u(S),\ \xi(T_u)>g_u(S+T_u)\Big\}\\\nonumber
&=\frac{1}{\sqrt{2\pi}}\int_{g_u(S)}^{\infty}\mathbb{P}\Big\{\mathcal{N}>\frac{g_u(S+T_u)-xr_{\xi}(T_u)}{\sqrt{1-r_{\xi}^2(T_u)}}\Big\}e^{-\frac{x^2}{2}}dx\\\label{fract}
&=\frac{1}{\sqrt{2\pi}g_u(S)}e^{-\frac{g_{u}^{2}(S)}{2}}\int_{0}^{\infty}
\mathbb{P}\Big\{\mathcal{N}>\frac{g_u(S+T_u)-\big(g_u(S)+x/g_u(S)\big)r_{\xi}(T_u)}{\sqrt{1-r_{\xi}^2(T_u)}}\Big\}
e^{{-\frac{x^2}{2g_{u}^{2}(S)}}-x}dx,
\end{align}
where $r_{\xi}(t)=e^{-Ct^2}$ is the correlation function of stationary Gaussian process $\xi(t)$. Furthermore, the fraction part within the probability in (\ref{fract}) can be rewritten as
\BQN\label{fraction1}
\frac{g_u(S)\big(\sigma(S)-\sigma(S+T_u)r_{\xi}(T_u)\big)+c\big(\tilde{\delta}(S+T_u)-\tilde{\delta}(S)\big)-\sigma(S+T_u)r_{\xi}(T_u)x/g_u(S)}
{\sigma(S+T_u)\sqrt{1-r_{\xi}^2(T_u)}}.
\EQN
With Taylor formula, simple calculations indicate \eqref{fraction1} tends to $-\infty$ as $u\to\infty$. Consequently, (\ref{fract}) is asymptotically equivalent to $\Psi\LT(g_u(S)\RT)$.\\
Next, thanks to the Bulinskaya's theorem, see Theorem E.4 in \cite{piterbarg1996asymptotic}, the probability of contingence of any level $u$ by the process $\hat{\xi}(\cdot)$ is equal to zero. Therefore, the event in (\ref{negli-prob}) implies that the number of crossings of the level $g_u(S)$ by the process $\hat{\xi}(\cdot)$ is greater than one. Denote by $N_{g_u(S)}[0,T_u]$ the number of crossings, by $p_{st}(\cdot,\cdot,\cdot,\cdot)$ the distribution density of the vector $(\hat{\xi}(s),\hat{\xi}(t),\hat{\xi}'(s),\hat{\xi}'(t))$ and by $p_{st}(\cdot,\cdot,\cdot)$ the distribution density of the vector $(\hat{\xi}(s),\hat{\xi}(t),\hat{\xi}'(t))$ . Then, appealing to the Theorem E.2 of \cite{piterbarg1996asymptotic} and using the symmetry of a normal density, we get a series of upper bounds for (\ref{negli-prob}),
\BQN
(\ref{negli-prob})&\leq&\pk{N_{g_u(S)}[0,T_u]\geq2}\nonumber\\\nonumber
&\leq& \mathbb{E}\big[N_{g_u(S)}[0,T_u]\big(N_{g_u(S)}[0,T_u]-1\big)\big]\\\nonumber
&=&\int_{0}^{T_u}\int_{0}^{T_u}\int_{-\infty}^{\infty}\int_{-\infty}^{\infty}\abs{xy}p_{st}(g_u(S),g_u(S),x,y)dxdydsdt\\\nonumber
&\leq&\int_{0}^{T_u}\int_{0}^{T_u}\int_{-\infty}^{\infty}\int_{-\infty}^{\infty}\frac{x^2+y^2}{2}p_{st}(g_u(S),g_u(S),x,y)dxdydsdt\\\nonumber
&\leq&\int_{0}^{T_u}\int_{0}^{T_u}\int_{-\infty}^{\infty}y^2 p_{st}(g_u(S),g_u(S),y)dydsdt\\\nonumber
&=&\int_{0}^{T_u}\int_{0}^{T_u}\varphi_{st}\big(g_u(S),g_u(S)\big)\int_{-\infty}^{\infty}y^2 \varphi_{st}\big(y|g_u(S),g_u(S)\big)dydsdt\\\label{negli-upp}
&=&\int_{0}^{T_u}\int_{0}^{T_u}\varphi_{st}\big(g_u(S),g_u(S)\big)\big(m_{u}^2(s,t)+\sigma^2_{u}(s,t)\big)dsdt,
\EQN
where $\varphi_{st}(\cdot,\cdot)$ represents the distribution density of $(\hat{\xi}(s),\hat{\xi}(t))$. Specifically,
\begin{align}
&\varphi_{st}\big(g_u(S),g_u(S)\big) \nonumber\\
&\quad=\frac{\sigma_{X_u}^{-1}(S+s)\sigma_{X_u}^{-1}(S+t)}{2\pi\sqrt{1-r_{\xi}^2(t-s)}}
\exp{\LT(-\frac{g_u^2(S)}{2}\Big[\frac{\LT(\sigma_{X_u}^{-1}(S+t)-r_\xi(t-s)\sigma_{X_u}^{-1}(S+s)\RT)^2}{1-r_{\xi}^2(t-s)}+\sigma_{X_u}^{-2}(S+s)\Big]\RT)}\label{negli-exp}.
\end{align}
The expected value $m_u(s,t)$ and the variance $\sigma^2_{u}(s,t)$ of the conditional density $\varphi_{st}(y|g_u(S),g_u(S))$ of random variable $\hat{\xi}'(t)$ given $\hat{\xi}(s)=\hat{\xi}(t)=g_u(S)$, according to Lemma \ref{lema-3gauss}, are equal to
\BQNY
m_u(s,t)=g_{u}(S)\LT(\frac{\sigma_{X_u}'(S+t)}{\sigma_{X_u}(S+t)}+\frac{r_{\xi}'(t-s)\big(\sigma_{X_u}(S+t)-r_\xi(t-s)\sigma_{X_u}(S+s)\big)}
{\sigma_{X_u}(S+s)(1-r_{\xi}^2(t-s))}\RT)
\EQNY
and
\BQNY
\sigma^2_u(s,t)=\sigma_{X_u}^2(S+t)\LT(2C-\frac{r_{\xi}'(t-s)^2}{1-r_{\xi}^2(t-s)}\RT).
\EQNY
Applying the Taylor formula again, after some technical calculations we have for $s,t\in[0,T_u]$
\BQNY
g_{u}^{-1}(S)\frac{m_u(s,t)}{\abs{t-s}}\to -C\qquad \textmd{and}\qquad \frac{\sigma^2_u(s,t)}{\abs{t-s}^2}\to 2C^2
\EQNY
uniformly as $u\to\infty$. Similar calculations show that the quantity in the square brackets of (\ref{negli-exp}) converges to $(2C)^{-1}\sigma'(S)^2/\sigma^2(S)+1$ uniformly as $u\to\infty$. Substituting these asymptotic results back into (\ref{negli-upp}) yields an upper bound in the form
\BQNY
Constant* T_{u}^{2}g_{u}^{2}(S)\exp\LT(-\frac{g_{u}^{2}(S)}{2}\LT(1+\frac{\sigma'(S)^2}{2C\sigma^2(S)}\RT)\RT),
\EQNY
which is negligible with respect to $\Psi\LT(g_u(S)\RT)$ as $u\to\infty$. \\
In summary, for all three different cases of $T_u$, we have the lower bound \eqref{lower-bound}. This together with upper bound (\ref{upp-bound}) completes the proof.\QED

\begin{remark}
The attentive reader may have found that the method of moments in Case III can be also applied to Case I and Case II. However, the method used in Case I and Case II, as the authors have tried before, failed to solve Case III. The method of moments, also known as Rice method, has long been used to estimate the distribution of the maximum of a random process (See the monograph \cite{leadbetter2012extremes} and recent paper \cite{kobelkov2005ruin}).
\end{remark}

\prooftheo{pan3}
From the definition of ruin time $\tau(u)$ in \eqref{eq:eta2}, we know
\BQNY
\pk{\tau(u)\leq S+T_u}=\mathcal{P}_S(u,T_u).
\EQNY
Then, by \eqref{psu-ano-form}, for any $x>0$
\BQNY
\pk{u^2(S+T_u-\tau(u))> x|\tau(u)<S+T_u}
=\frac{\pk{\sup_{t\in[0,S_x(u)]}\inf_{s\in[t,t+T_{u}]}\widetilde{X}_u(s)>g_u(S_x(u))}}{\pk{\sup_{t\in[0,S]}\inf_{s\in[t,t+T_{u}]}{X}_u(s)>g_u(S)}},
\EQNY
with $S_x(u):=S- x u^{-2}$ and $\widetilde{X}_u(t):=\frac{Y(t)}{\sigma(t)}\frac{g_u(S_x(u))}{g_u(t)}$. We need to find an exact asymptotic for the numerator.
First, as in the proof of Theorem \ref{pan1} case III, just replacing $S$ by $S_x(u)$, we have
\BQNY
\mathcal{P}_{S_x(u)}(u,T_{u}) \geq \pk{\inf_{s\in[S_x(u),S_x(u)+T_{u}]}\widetilde{X}_u(s)>g_u(S_x(u))} \geq \Psi\big(g_u(S_x(u))\big)\ooo
\EQNY
as $u\to \IF$. Next, appealing to the upper bound given in the proof of Theorem 2.4 in \cite{dbicki2015gaussian}, we get
\BQNY
\mathcal{P}_{S_x(u)}(u,T_{u}) \leq \pk{\sup_{t\in[0,S_x(u)]}\widetilde{X}_u(t)>g_u(S_x(u))} = \Psi\big(g_u(S_x(u))\big)\ooo
\EQNY
as $u\to \IF$.
 Therefore,
\begin{eqnarray*}
\pk{u^2(S+T_u-\tau(u))> x|\tau(u)<S+T_u}&\sim&\frac{\Psi(g_u(S_x(u)))}{\Psi(g_u(S))}\\
&\sim&\exp\left(\frac{g_u^2(S)-g_u^2(S_x(u))}{2} \right), \quad u \to \IF.
\end{eqnarray*}
Some standard algebra yields 
\BQNY
g_u^2(S)-g_u^2(S_x(u)) \rightarrow -\frac{2\sigma'(S)}{\sigma^3(S)}x, \quad \textmd{as} \ u\to\IF.
\EQNY
In other words,
\begin{eqnarray*}
\lim_{u\to\infty}
\pk{u^2(S+T_u-\tau(u))> x|\tau(u)<S+T_u}
=\exp\left(  - \frac{\sigma'(S)}{\sigma^3(S)}x \right),
\end{eqnarray*}
which completes the proof.                    \QED

%

\BEL\label{lema-3gauss}
Let $(X,Y,Z)$ be a centered Gaussian vector with values in $\mathcal{R}^3$, the conditional distribution of $Z$ given $X=x$ and $Y=y$ is a Gaussian random variable with expected value
\BQNY
\mathbb{E}[Z|X=x,Y=y]=(x,y)Q^{-1}\textbf{b}
\EQNY
and variance
\BQNY
\mathbb{V}ar(Z|X=x,Y=y)=Var(Z)-\textbf{b}^{T}Q^{-1}\textbf{b},
\EQNY
where $Q$ is the covariance matrix of random variables $X$ and $Y$, $\textbf{b}=(\mathbb{C}ov(X,Z),\mathbb{C}ov(Y,Z))^{T}$.
\EEL
\prooflem{lema-3gauss}
 Decomposing $Z$ as sum like
\BQNY
Z=\alpha X+\beta Y+ \Gamma
\EQNY
such that $\Gamma$ is independent of both $X$ and $Y$. This yields $(\alpha,\beta)^T=Q^{-1}\textbf{b}$.     \QED\\

{\bf Acknowledgment} X. Peng was partially supported by National Natural Science Foundation of China (71501025), the Fundamental Research Funds for the Central Universities (ZYGX2015J102), and Chinese Scholarship Council. This work was completed when the authors were visiting
the University of Lausanne. Both the authors kindly acknowledge partial support by the SNSF Grant 200021-166274.

\bibliographystyle{plain}
\bibliography{MMM}
\end{document}